\documentclass[a4paper,twoside,11pt,reqno]{amsart}

\usepackage[headings]{fullpage}
\usepackage{lmodern}
\usepackage{microtype}
\usepackage{amsmath}
\usepackage{amsthm}
\usepackage{amssymb}
\usepackage{xcolor}
\usepackage[hypertexnames=false]{hyperref}
\usepackage{thmtools}
\hypersetup{colorlinks=true,linkcolor=red,citecolor=black,urlcolor=black}
\usepackage{mathrsfs}
\usepackage{enumitem}
\usepackage{tikz-cd}
\usepackage{aliascnt}
\usepackage[capitalize]{cleveref}

\theoremstyle{definition}
\newtheorem{defi}{Definition}[section]

\newtheorem{rem}[defi]{Remark}
\theoremstyle{plain}
\newtheorem*{thm*}{Theorem}
\newtheorem*{cor*}{Corollary}
\newtheorem{thm}[defi]{Theorem}
\newtheorem{lem}[defi]{Lemma}
\newtheorem{cor}[defi]{Corollary}
\newtheorem{prop}[defi]{Proposition}

\newtheorem*{theoremA}{Theorem A}
\newtheorem*{theoremB}{Theorem B}
\newtheorem*{theoremC}{Theorem C}
\newtheorem*{theoremD}{Theorem D}

\newcommand{\Alb}{\mathrm{Alb}}
\newcommand{\C}{\mathbb{C}}
\newcommand{\Q}{\mathbb{Q}}
\newcommand{\Z}{\mathbb{Z}}
\newcommand{\F}{\mathbb{F}}
\newcommand{\X}{\mathcal{X}}
\newcommand{\Y}{\mathcal{Y}}
\newcommand{\nr}{\mathrm{nr}}
\newcommand{\Cl}{\mathrm{Cl}}
\newcommand{\p}{\mathfrak{p}}
\newcommand{\ab}{\mathrm{ab}}
\newcommand{\tors}{\mathrm{tors}}

\newcommand{\geo}{\mathrm{geo}}
\newcommand{\lisim}{\overset{\sim}{\longrightarrow}}
\newcommand{\iso}{\simeq}
\renewcommand{\o}{\mathcal{O}}
\renewcommand{\H}{\mathrm{H}}
\renewcommand{\phi}{\varphi}
\renewcommand{\l}{\ell}
\DeclareMathOperator{\spec}{Spec}
\DeclareMathOperator{\im}{im}
\DeclareMathOperator{\Gal}{Gal}
\DeclareMathOperator{\Frob}{Frob}
\DeclareMathOperator{\Aut}{Aut}
\newcommand{\bp}{\begin{proof}}
\newcommand{\ep}{\end{proof}}
\newcommand{\be}{\begin{enumerate}}
\newcommand{\ee}{\end{enumerate}}

\newcommand\restr[2]{{
		\left.\kern-\nulldelimiterspace 
		#1 
		\vphantom{\big|} 
		\right|_{#2} 
}}

\begin{document}
\title{Galois coinvariants and fundamental groups of arithmetic schemes}
\author{Benjamin Steklov}
\address{Benjamin Steklov, Institut f\"ur Mathematik, Goethe--Universit\"at Frankfurt, Robert-Mayer-Str.~6--8, 60325~Frankfurt am Main, Germany}
\email{steklov@math.uni-frankfurt.de}
\date{\today} 

\begin{abstract}
We investigate the structure of étale fundamental groups of arithmetic schemes over rings of integers of number fields. We show that the image of the geometric étale fundamental group is FAb, while the kernel of the map to the fundamental group of the base need not be FAb or finitely generated. Our results rely on finiteness theorems for Galois coinvariants of geometric étale fundamental groups and étale cohomology over finitely generated fields. These apply to Galois extensions $L/k$ for which $\Gal(L/k)^\ab$ has finite exponent. As a further application, we establish finiteness of prime-to-characteristic $L$-rational torsion on abelian varieties over $k$.
\end{abstract}
\maketitle
\tableofcontents

\section{Introduction}

Let $F$ be a number field with ring of integers $\o_F$.
Bost--Charles \cite{bostQuasiprojectiveFormalanalyticArithmetic2026}
investigated when the homomorphism
\[
x_*\colon \pi_1(\o_F)\to\pi_1(\X)
\]
induced by a section $x\in\X(\o_F)$ of an arithmetic
surface $\X/\o_F$ has image of finite index. They also asked whether there
exists a regular arithmetic surface $\X/\Z$ admitting
a section and having infinite étale fundamental group.
This question was recently answered affirmatively by
Chenal \cite{chenalProperArithmeticSurface2026}
and Keller--Pagano
\cite{kellerGolodShafarevichArithmeticSurfaces2026}.

In this paper, we investigate finiteness properties of
fundamental groups of arithmetic schemes in a more general
setting, without assuming the existence of a section.

Let $\X$ be a normal scheme, flat, surjective and
of finite type over $\o_F$, with geometrically
connected generic fiber $\X_F$.
Associated with $\X/\o_F$ is a complex of
fundamental groups
\[
\pi_1(\X_{\overline F})\to
\pi_1(\X)\to
\pi_1(\o_F)\to 1,
\]
which need not be exact in the middle,
see \cref{diagram}. We therefore distinguish the two subgroups
\[
\pi_1^{\geo}(\X/\o_F)
:=
\im\bigl(\pi_1(\X_{\overline F})\to\pi_1(\X)\bigr)
\subset
\ker\bigl(\pi_1(\X)\to\pi_1(\o_F)\bigr)
=:K(\X/\o_F).
\]
The inclusion is an equality whenever $\X/\o_F$ admits a section,
see \cref{lem:sectionimpliesexactness}.
In this case, the finite-index problem of Bost--Charles
amounts to asking whether $K(\X/\o_F)$ is finite.
The examples of Chenal and Keller--Pagano show that
both groups can be infinite, since
\[
\pi_1^{\geo}(\X/\Z)=K(\X/\Z)=\pi_1(\X)
\]
for their arithmetic surfaces. It is therefore natural to investigate the structure
of these two groups separately.

Recall that a profinite group is called \textit{FAb}
if every open subgroup has finite abelianization.
By class field theory, $\pi_1(\o_F)$ is FAb,
and a theorem of Schmidt
\cite[Corollary 3.3]{schmidtTameCoveringsArithmetic2002}
shows that $\pi_1(\X)$ is FAb as well.
However, the FAb property need not pass to closed
subgroups, so these results do not directly apply to
$\pi_1^{\geo}(\X/\o_F)$ or $K(\X/\o_F)$.

The group $\pi_1^{\geo}(\X/\o_F)$ is finitely
generated, being a quotient of
$\pi_1(\X_{\overline F})$.
Our first result shows that it is also FAb.
\begin{theoremA}[See \cref{thm:arithmeticfab}]
	Let $\X$ be a normal scheme, flat, surjective
	and of finite type over $\o_F$, with
	geometrically connected generic fiber $\X_F$.
	Then $\pi_1^{\geo}(\X/\o_F)$ is FAb.
\end{theoremA}
Theorem A also provides an alternative proof of
Schmidt's result mentioned above, see
\cref{cor:schmidt}. By contrast, the larger group $K(\X/\o_F)$
can have a fundamentally different structure.

\begin{theoremB}[See \cref{thm:infinitekernel}]
	There exists an imaginary quadratic number field
	$F$ and infinitely many regular proper arithmetic
	surfaces $\X/\o_F$ such that the pro-$2$ completion
	$K(\X/\o_F)^{(2)}$
	is not finitely generated.
\end{theoremB}

In particular, $K(\X/\o_F)$ need not be FAb
or finitely generated. By Theorem A, the surfaces
in Theorem B cannot admit a section.

The proof of Theorem B uses a theorem of
Hajir--Maire--Ramakrishna
\cite[Theorem 4]{hajirCuttingTowersNumber2021}.
Theorem A, on the other hand, is a consequence
of a more general finiteness result for Galois
coinvariants of geometric étale fundamental groups
over finitely generated fields, which we now describe.

\subsection{Galois coinvariants}

Let $X$ be a geometrically connected variety
over a field $k$, with separable closure $\bar{k}$
and geometric point $\bar{x}\in X_{\bar{k}}$.
The fundamental exact sequence
\[
1\longrightarrow\pi_1(X_{\bar{k}},\bar{x})
\longrightarrow\pi_1(X,\bar{x})
\longrightarrow\Gal_k\longrightarrow 1
\]
induces a canonical $\Gal_k$-action on
$\pi_1(X_{\bar{k}})^\ab$. A section
$s\colon\Gal_k\to\pi_1(X,\bar{x})$ of this sequence induces a $\Gal_k$-action on the
group $\pi_1(X_{\bar{k}},\bar{x})$ via
\[
\sigma.\pi:=s(\sigma)\pi s(\sigma)^{-1},
\qquad
\sigma\in\Gal_k,\quad
\pi\in\pi_1(X_{\bar{k}},\bar{x}).
\]
For any separable algebraic extension
$L/k$, we may restrict
these actions to $\Gal_L\subset\Gal_k$
and consider the coinvariants
\[
\left[\pi_1(X_{\bar{k}})^\ab\right]_{\Gal_L}
\quad\text{and}\quad
\left[\pi_1(X_{\bar{k}},\bar{x})\right]_{\Gal_L},
\]
where only the latter require a choice of section.

We obtain the following finiteness result.
\begin{theoremC}[See \cref{thm:fabnessoffundamentalgroups}]
	Let $X$ be a geometrically connected and
	geometrically unibranch variety over a finitely
	generated field $k$ of characteristic $p\geq 0$.
	Let $L/k$ be a (possibly infinite) Galois
	extension such that $\Gal(L/k)^\ab$ has
	finite exponent.
	\be
	\item The group of coinvariants
	\[
	\left[\pi_1(X_{\bar{k}})^{\ab,(p')}\right]_{\Gal_L}
	\]
	is finite.
	
	\item Let $\bar{x}\in X_{\bar{k}}$ be a geometric
	point and $s\colon\Gal_k\to\pi_1(X,\bar{x})$
	a section of the fundamental exact sequence.
	Endow $\pi_1(X_{\bar{k}},\bar{x})$ with the
	$\Gal_k$-action induced by $s$.
	For any open subgroup
	\[
	H\subset
	\left[\pi_1(X_{\bar{k}},\bar{x})\right]_{\Gal_L},
	\]
	the group $H^{\ab,(p')}$ is finite.
	In particular, if $p=0$, the group
	\[
	\left[\pi_1(X_{\bar{k}},\bar{x})\right]_{\Gal_L}
	\]
	is FAb.
	\ee
\end{theoremC}

We briefly explain how Theorem C implies
Theorem A. By a theorem of Rumely and
Moret-Bailly
\cite[Théorème 0.7]{moret-baillyProblemesSkolemChamps2001},
we may reduce to the case where
$\X/\o_F$ admits a section $x\in\X(\o_F)$.
Let $F^\nr/F$ denote the maximal extension
unramified at all finite places, so that
$\Gal(F^\nr/F)=\pi_1(\o_F)$.
The point $x$ endows
$\pi_1(\X_{\overline F},x)$ with a $\Gal_F$-action,
and there is a natural surjection
\[
\left[
\pi_1(\X_{\overline F},x)
\right]_{\Gal_{F^\nr}}
\twoheadrightarrow
\pi_1^{\geo}(\X/\o_F).
\]
By class field theory, the abelianization
$\Gal(F^\nr/F)^\ab$ is isomorphic to the narrow
class group of $F$, which is finite.
Thus, Theorem C applies.

The proof of Theorem C follows a strategy
similar to that of Katz--Lang
\cite{katzFinitenessTheoremsGeometric1981},
reducing the problem to semi-abelian varieties.
A key ingredient is a finiteness criterion for
coinvariants of systems of pure Galois
representations of nonzero weight, see
\cref{prop:coinvariantsofpure}.
This criterion also yields a uniform finiteness
result for the étale cohomology of smooth,
proper varieties.

\begin{theoremD}[See \cref{thm:cohomology}]
	Let $X$ be a smooth, proper variety over a
	finitely generated field $k$ of characteristic
	$p\geq 0$. Fix integers $i,j$ with $i\neq 2j$
	and $e\geq 1$.
	For every (possibly infinite) Galois extension
	$L/k$ with $e\Gal(L/k)^\ab=0$, the group
	of coinvariants
	\[
	\left[
	\H^i(X_{\bar{k}},\widehat{\Z}^{(p')}(j))
	\right]_{\Gal_L}
	\]
	is finite. Moreover, its order divides
	an integer $D=D(X/k,i,j,e)>0$
	independent of $L$.
\end{theoremD}

As an application, we obtain the following
uniform finiteness result for torsion points
on abelian varieties.

\begin{cor*}[See \cref{cor:torsionpoints}]
	Let $A$ be an abelian variety over a finitely
	generated field $k$ of characteristic $p\geq 0$
	and let $e\geq 1$.
	For every (possibly infinite) Galois extension
	$L/k$ with $e\Gal(L/k)^\ab=0$, the group
	\[
	A(L)_{\tors}^{(p')}
	\]
	is finite and its order divides an integer
	$D=D(A/k,e)>0$ independent of $L$.
\end{cor*}

The bound $D$ can be expressed in terms of
characteristic polynomials of Frobenii acting
on $\l$-adic cohomology. In the case of elliptic
curves, it can be made explicit, see
\cref{cor:torsionelliptic}.

\subsection{Related work}

The finiteness theorems of Katz--Lang
\cite{katzFinitenessTheoremsGeometric1981}
are closely related to Theorem C.
In fact, their strategy can be adapted
to prove Theorem C in the special case
$L=k$. Extending this approach to infinite
Galois extensions $L/k$, however, requires
new arguments. Nonabelian Galois coinvariants
for $L=k$ were also studied by
Frey--Kani--Völklein
\cite{freyCurvesInfiniteKrational1999a}
under the name of $k$-\textit{rational
	geometric fundamental groups}.

The statement and proof of Theorem D
have parallels with the work of Rössler--Szamuely
\cite{rosslerCohomologyTorsionCycles2019}
and Lombardo--Szamuely
\cite{lombardoCohomologyVarietiesMaximal2026}.
Their results concern the maximal cyclotomic
and maximal Kummer extensions of number fields,
respectively. Both extensions have Galois groups
whose abelianizations are of infinite exponent,
whereas Theorem D applies to Galois extensions
$L/k$ such that $\Gal(L/k)^\ab$ has finite exponent.
Thus, the results are complementary and concern
different classes of Galois extensions.

Finally, let $A/F$ be an abelian variety
with Néron model $\mathcal A/\o_F$.
Since $\pi_1(A_{\overline F})$ is abelian,
the group $\pi_1^{\geo}(\mathcal A/\o_F)$
is abelian and hence finite by Theorem A.
This recovers a result recently obtained
by Lu \cite[Theorem 1.1]{luFinitenessEtaleFundamental2026}
using Faltings heights.
We refer to
\cite{luFinitenessEtaleFundamental2026}
for a discussion of uniformity and further
results on the order of
$\pi_1^{\geo}(\mathcal A/\o_F)$.

\subsection*{Acknowledgments} I would like to thank my PhD advisor Jakob Stix for bringing the question posed in \cite{bostQuasiprojectiveFormalanalyticArithmetic2026} to my attention and many helpful discussions. I further thank Magnus Carlsson, Simon Kaib, Immanuel Klevesath, Alexander Schmidt, Ruth Wild and Nils Witt for helpful comments and conversations related to this work. 

The present article is part of the PhD thesis of the author. This project was funded by Deutsche Forschungsgemeinschaft (DFG, German Research Foundation) through the Collaborative Research
Centre TRR 326 \textit{Geometry and Arithmetic of Uniformized Structures} - Project-ID 444845124.
\subsection*{AI Disclosure} ChatGPT was used for finding references and learning about existing results. During the drafting process ChatGPT was used for stylistic improvements of the introduction. The article was written by the author and the proofs are due to the author. The author takes full responsibility for the correctness of the content.
\subsection*{Notation.}
 A \textbf{variety} over a field $k$ is a finite type scheme $X/k$. A separable closure of $k$ is denoted by $\bar{k}$. The letter $p$ always denotes either a rational prime or $0$. For a profinite group $G$, we  write $G^{(p')}$ (resp. $G^{(p)}$) for the prime-to-$p$ completion (resp. pro-$p$-completion), where $G^{(0')}:=G$.  We write $\pi_1(X,\bar{x})$ for the \'etale fundamental group of a scheme $X$ based at a geometric point $\bar{x}$ or simply $\pi_1(X)$, when the choice of $\bar{x}$ is irrelevant. For a ring $R$ with connected spectrum, we write $\pi_1(R):=\pi_1(\spec(R))$.
 \section{Galois representations}
\subsection{An abstract finiteness criterion for coinvariants}
We begin by recalling the notion of a pure $\l$-adic Galois representation of a finite field $\F$, following Deligne's convention that $\Frob_\F$ denotes the geometric Frobenius.
\begin{defi}
	Let $\F$ be a finite field of characteristic $p>0$ and let $\l \neq p$ be a prime. An $\l$-adic representation $V$ of $\Gal_\F$ is called \textbf{pure} of weight $w \in \Z$ if all roots $\alpha \in \overline{\Q}_\l$ of the characteristic polynomial 
	\[ \det(X- \Frob_\F\mid V ) \in \Q_\l[X]\]
	are algebraic over $\Q$ and have absolute value $(\# \F)^{w/2}$ under every embedding $\iota \colon \Q(\alpha)\to \C$. We emphasize that we do not require the characteristic polynomial to be rational.
\end{defi}
For a non-zero rational number $a$ and $\Sigma$ a set of rational primes, we denote by
\[ a_\Sigma:= \prod_{\l \in \Sigma}^{}\l^{v_\l(a)} \in \Q_{>0}\]
the $\Sigma$-part of $a$.
\begin{prop}\label{prop:coinvariantsofpure}
Let $k$ be a finitely generated field of characteristic $p\geq 0$. Let $\Sigma$ be a set of rational primes not containing $p$ and let $N\geq 1$.	For each  $\l \in \Sigma$, let $\Lambda_\l$ be a free $\Z_\l$-module of rank $N$ equipped with a continuous $\Gal_k$-action and set $V_\l:=\Lambda_\l \otimes \Q_\l$. Assume there exists a finite covering 
\[\bigcup_{r=1}^{s}\Sigma_r = \Sigma\]
such that for each $1\leq r\leq s$ there is an integral scheme $T_r$ of finite type over $\Z$ with function field $k$ and a closed point $t_r \in T_r$ satisfying the following conditions.
\begin{itemize}
	\item For every  $\l  \in \Sigma_r$ the representation $V_\l$ factors as a representation of $\pi_1(T_r[1/\l])$.
	\item  The residue characteristic of $t_r$ does not belong to $\Sigma_r$ and the representation
	\[\restr{V_\l}{\Gal_{\kappa(t_r)}}\]
is pure of non zero-weight for every $\l \in \Sigma_r$. Moreover the characteristic polynomial
	\[\det(X-\Frob_{t_r}\mid V_\l) \in \Q[X]\]
	is rational and independent of the choice $\l \in \Sigma_r$.
\end{itemize}
For an integer $e\geq 1$ and $1\leq n \leq N$ define the polynomials
\[P_{n,r}:=\det(X- \Frob_{t_r}^e \mid \bigwedge^n V_\l ) \in \Q[X],\qquad \l \in \Sigma_r.  \]
Then, for any (possibly infinite) Galois extension $L/k$ such that $e\Gal(L/k)^\ab=0$, the group of coinvariants
	\[\left[\prod_{\l \in \Sigma} \Lambda_\l \right]_{\Gal_L}  \]
	is finite and 
	divides the integer
	\[D(e):= \left(\prod_{r=1}^{s} \prod_{n=1}^{N}P_{n,r}(1)_{\Sigma_r}\right)^N>0.\]
\end{prop}
\bp
It suffices to treat the case of $s=1$. To simplify notation, we write $t=t_1 \in T_1$.

Step 1: We claim that the coinvariants $\left[\Lambda_\l\right]_{\Gal_L}$ are finite for all $\l \in \Sigma $.
Let $\l \in \Sigma$ and set $W_\l:= \left[V_\l\right]_{\Gal_L}$. It suffices to prove that $W_\l=0$. Suppose that $W_\l\neq 0$. As a $\Gal_{\kappa(t)}$-representation $V_\l$ is pure of weight $w\neq 0$,  therefore the quotient $W_\l$ is too. The determinant representation 
\[\Gal_k\to \Aut(\det(W_\l))=\Q_\l^\times\] factors through $\Gal(L/k)^\ab$. The image has finite exponent, hence is finite. Therefore 
\[ 0= \mathrm{weight}(\det(\restr{W_\l}{\Gal_{\kappa(t)}}))=\mathrm{weight}(\restr{W_\l}{\Gal_{\kappa(t)}}) \dim W_\l=w \dim W_\l \neq 0,\]
which is a contradiction.

Step 2:  By assumption, the polynomial
\[ \det(X-\Frob_{t}\mid V_\l) \in \Q[X]\]
is independent of the choice $\l \in \Sigma$, has degree $N\geq 1$ and all roots have absolute value 
\[ (\#\kappa(t))^{w/2},\]
where $w\neq 0$ is the weight of $\restr{V_\l}{\Gal_{\kappa(t)}}$.   As a consequence, for any $1\leq n\leq N$, the polynomials
\[P_n:=P_{n,1}=\det(X-\Frob_{t}^e \mid \bigwedge^nV_\l)\in \Q[X] \]
are rational, independent of the choice $\l \in \Sigma$ and all roots have absolute value 
\[(\# \kappa(t))^{nw/2} \neq 1.\] In particular, $P_n(1)$ is a non-zero rational number. More precisely, we have 
\[P_{n} \in \Q[X] \cap \bigcap_{\l \in \Sigma}\Z_\l[X]= \bigcap_{\l \in \Sigma} \Z_{(\l)}[X],\]
which implies that $P_n(1)_\Sigma$ is a non-zero integer.   
For $\l \in \Sigma$, we set $M_\l:= \left[\Lambda_\l\right]_{\Gal_L}$.  Because of Step 1, $M_\l$ is a finite abelian $\l$-group, hence we have a decomposition (as abelian groups)
 \[M_\l \iso \bigoplus_{j=1}^{d}(\Z/\l^{a_j}\Z),\]
where $d$ and $a_j$ are non-negative integers such that
 \[a_1\geq a_2 \geq \cdots \geq a_d \geq 0. \]
 The integers $a_j$ and $d$ may depend on $\l$. For $d$ we have the bound
 \[d=\dim_{\F_\l} M_\l \otimes \F_\l \leq \dim_{\F_\l} \Lambda_\l \otimes \F_\l =N.\]
 For $1\leq n\leq d$, consider the $\Gal_k$-equivariant quotient 
 \[M_\l \twoheadrightarrow M_\l(n):=\l^{a_{n+1}}M/\l^{a_n}M, \quad x \longmapsto \l^{a_{n+1}}x.\]
 Let $b_n:=a_{n}-a_{n+1}$, where $a_{d+1}:=0$. Then $M_\l(n) $ is a free $(\Z/\l^{b_n})$-module of rank $n$. We obtain the following $\Gal_k$-equivariant surjections
 \[\bigwedge^n \Lambda_\l \twoheadrightarrow \bigwedge^n M_\l \twoheadrightarrow  \bigwedge^{n}M_\l(n)=\det(M_\l(n)).\]
The action of $\Gal_k$ on $\det(M_\l(n))$ factors through the quotient $\Gal(L/k)^\ab$, which satisfies $e\Gal(L/k)^\ab=0$. Hence $\Frob_t^e$ acts trivially on $\det(M_\l(n))$. Thus we obtain $ \det(M_\l(n))$ as a quotient of the coinvariants
\[ \left[\bigwedge^n\Lambda_\l\right]_{\Frob_{t}^e} \twoheadrightarrow\det(M_\l(n)), \]
which implies
 \[b_n = v_\l \left(\# \det(M_\l(n))\right) \leq  v_\l \left(\#  \left[\bigwedge^n\Lambda_\l\right]_{\Frob_{t}^e}\right) =v_\l(P_n(1)).  \] 
We have
\[a_1 = \sum_{n=1}^{d}b_n\leq  \sum_{n=1}^{d} v_\l(P_n(1))\leq \sum_{n=1}^{N} v_\l(P_n(1))\]
and 
\[v_\l(\# M_\l) =\sum_{n=1}^{d}a_n\leq d a_1 \leq N \sum_{n=1}^{N} v_\l(P_n(1)). \]
Therefore the number
 \[
 \# \left[\prod_{\l \in \Sigma} \Lambda_\l\right]_{\Gal_L} \] 
 divides the integer
 \[D(e):= \left(\prod_{n=1}^NP_{n}(1)_\Sigma\right)^N> 0 .\qedhere\]

\ep 

\subsection{\'Etale cohomology and torsion points}
\begin{thm}\label{thm:cohomology}
	Let $X$ be a smooth, proper variety over a finitely generated field $k$ of characteristic $p\geq 0$. Fix integers $i,j$ with $i\neq 2j$ and $e\geq 1$. For every (possibly infinite) Galois extension $L/k$ with $e\Gal(L/k)^\ab =0$ the group of coinvariants 
\[ \left[\H^i(X_{\bar{k}},\widehat{\Z}^{(p')}(j))\right]_{\Gal_L}  \]
is finite. Moreover, its order divides an integer $D=D(X/k,i,j,e)>0$ independent of $L$.
\end{thm}
\bp 
By a theorem of Gabber \cite{gabberTorsionDansCohomologie1983} and de Jongs alterations \cite{dejongSmoothnessSemistabilityAlterations1996}, the \'etale cohomology groups
\[  \H^i(X_{\bar{k}},\Z_\l(j)) \] are torsion-free for almost all $\l\neq p$. Hence 
\[ \H^i(X_{\bar{k}},\widehat{\Z}^{(p')}(j))\twoheadrightarrow \prod_{\l\neq p} \H^i(X_{\bar{k}},\Z_\l(j))/ \text{torsion},\]
has finite kernel. Since applying $\Gal_L$-coinvariants is right-exact, it suffices to consider
\[ \left[\prod_{\l\neq p} \H^i(X_{\bar{k}},\Z_\l(j))/ \text{torsion}\right]_{\Gal_L}.\]
For $\l\neq p$, let
\[\Lambda_\l:=\H^i(X_{\bar{k}},\Z_\l(j))/\text{torsion}\]
and
\[N:= \mathrm{rank}_{\Z_\l} \Lambda_\l= \dim_{\Q_\l} \Lambda_\l\otimes \Q_\l = \dim_{\Q_\l} \H^i(X_{\bar{k}},\Q_\l).\]
By smooth proper base-change and \cite[Corollaire (3.3.9).]{deligneConjectureWeilII1980} the number $N$ is independent of the choice of $\l \neq p$. The case $N=0$ is trivial, so we may suppose that $N\geq 1$. It remains to argue that the assumptions of \cref{prop:coinvariantsofpure} apply to $\Lambda_\l$.
This is essentially a consequence of the main results of \cite{deligneConjectureWeilII1980}. To be precise, after spreading out we find a normal integral $T$ of finite type over $\Z$ with function field $k$ and a smooth proper morphism $\X/T$ with generic fiber $X/k$. Let $\Sigma$ be the set of all primes not equal to $p$. Then for any $\l$, the Galois representations 
\[V_\l:=\Lambda_\l\otimes \Q_\l:= \H^i(X_{\bar{k}},\Q_\l(j))=\Lambda_\l\otimes \Q_\l \] factor as a representations of $\pi_1(T[1/\l])$.

Case 1: $p>0$. We choose a closed point $t_1 \in T$ of residue characteristic $p$ and define $\Sigma_1:=\Sigma$.

Case 2: $p=0$. We choose closed points $t_1,t_2 \in T$ with distinct residue characteristics $\l_r:= \mathrm{char}(\kappa(t_r))$ and define 
\[\Sigma_r:= \Sigma \setminus \{\l_r\}\quad r=1,2.\]
In both cases we have $\Sigma = \bigcup_r \Sigma_r$. Let $\l \in \Sigma_r$. By smooth proper base change, we have 
\[ \restr{V_\l}{\Gal_{\kappa(t_r)}}=\H^i(\X_{\overline{\kappa(t_r)}},\Q_\l(j)).\]
By \cite[Corollaire (3.3.9).]{deligneConjectureWeilII1980} this representation is pure of weight $i-2j\neq 0$ with Frobenius having rational characteristic polynomial, independent of $\l\in \Sigma_r$. Hence \cref{prop:coinvariantsofpure} applies to $\Lambda_\l$.
\ep 
\begin{cor}\label{cor:torsionpoints}
	Let $A$ be an abelian variety over a finitely generated field $k$ of characteristic $p\geq 0$ and let $e\geq 1$. For every (possibly infinite) Galois extension $L/k$ with $e\Gal(L/k)^\ab=0$ the group
\[ A(L)_{\tors}^{(p')}\]
is finite and its order divides an integer $D=D(A/k,e)>0$ independent of $L$.
\end{cor}
	\bp 
The group $A(\bar{k})_{\tors}^{(p')}$ is Pontryagin dual to  $\H^1(A_{\bar{k}},\widehat{\Z}^{(p')})$. Hence $A(L)_{\tors}^{(p')}$ is Pontryagin dual to the coinvariants
\[\left[\H^1(A_{\bar{k}},\widehat{\Z}^{(p')})\right]_{\Gal_L}.\] Therefore \cref*{thm:cohomology} implies the claim.
\ep 
To obtain a concrete example for the bound $D$ we may consider elliptic curves.
\begin{cor}\label{cor:torsionelliptic}
	Let $E$ be an elliptic curve over a number field $F$ and  $\mathfrak{q}$ be a prime of good reduction with residue characteristic $\l_0$ and norm $q=N(\mathfrak{q})$. For any Galois extension $L/F$ with $e \Gal(L/F)^\ab =0$ we have
	\[  \# E(L)_{\tors}^{(\l_0)'} \mid  \left(\# E(\F_{q^e})(q^e-1)\right)^2.  \] 
\end{cor}
\bp 
Let $\Sigma$ be the set of all rational primes not containing $\l_0$. Applying \cref{prop:coinvariantsofpure} to 
\[\Lambda_\l:= \H^1(E_{\overline{F}},\Z_\l)\\\quad \text{for $\l \in \Sigma$}\]
implies that 
\[ \# E (L)^{(\l_0)'}_{\tors}= \# \left[ \H^1(E_{\overline{F}},\widehat{\Z}^{(\l_0)'})\right]_{\Gal_L}\]
divides the prime-to-$\l_0$ part of
\[ \left(P_1(1)P_2(1)\right)^2,\]
where 
\[P_n:= \det(X- \Frob_{\mathfrak{q}}^e \mid \bigwedge^n \H^1(E_{\overline{F}},\Z_\l) ) \quad \text{for $\l \neq \l_0$}.\]
We have 
\[P_{1}:=X^2-a_{q^e}X +q^e \text{ and  } P_{2}=X- q^e,\]
where 
\[ a_{q^e}:= q^e+1- \# E(\F_{q^e}). \]
Hence $\# E (L)_{\tors}^{(\l_0)'}$ divides (the prime-to-$\l_0$ part of) the integer
\[ \left(P_1(1)P_2(1)\right)^2=\left(\# E(\F_{q^e})(q^e-1)\right)^2.  \qedhere\]
\ep 
\begin{rem}
	Note that in \cref{cor:torsionelliptic} the role of $\l_0$ is not the same as the role of $p$ in \cref{cor:torsionpoints}. The choice of another prime of good reduction with residue characteristic not equal to $\l_0$ yields a bound on the order of the full group $E(L)_{\tors}$.
\end{rem}
\section{\'Etale fundamental groups}
\subsection{The Albanese torsor}
In this subsection, we recall the Albanese torsor of a variety, cf.\ \cite[Section 2]{wittenbergAlbaneseTorsorsElementary2008}.
\begin{defi}
	Let $X$ be variety over a field $k$. The \textbf{Albanese torsor} of $X/k$ is a triple $(\Alb^1_{X/k},\Alb^0_{X/k},u_{X/k})$, where  $\Alb^0_{X/k}$ is a semi-abelian variety over $k$, $\Alb^1_{X/k}$ is a torsor under $\Alb^0_{X/k}$ and 
	\[u_{X/k}\colon X\to \Alb^1_{X/k} \]
	is a $k$-morphism,  satisfying the following universal property. For any semi-abelian variety $A^0$ over $k$, any torsor $A^1$ under $A^0$ and any $k$-morphism $m\colon X\to A^1$, there exists a unique $k$-morphism $f^1\colon \Alb^1_{X/k}\to A^1$, such that $f^1\circ u_{X/k}=m$, and there exists a unique $k$-group homomorphism $f^0\colon \Alb^0_{X/k}\to A^0$, such that $f^1$ is $f^0$-equivariant.
\end{defi} 
We have the following theorem of Serre in the case where $k$ is algebraically closed, see \cite[Appendix A]{wittenbergAlbaneseTorsorsElementary2008} for the general case.
\begin{thm}[Serre]
	Let $X$ be a geometrically integral variety over a field $k$. Then the Albanese torsor of $X/k$ exists.
\end{thm}
Let $X/k$ be a geometrically integral variety and $K/k$ a field extension. The base change 
\[(u_{X/k})_K\colon X_K \longrightarrow \Alb^1_{X/k}\times_kK\]
implies the existence of unique $K$-morphisms 
\[ \theta^n_{K/k} \colon \Alb^n_{X_{K}/K} \longrightarrow \Alb^n_{X/k}\times_k K, \qquad n=0,1 \]
such that  $\theta^0_{K/k} $ is a group homomorphism, $\theta^1_{K/k}$ is $\theta^0_{K/k}$-equivariant and
\[\theta^1_{K/k} \circ u_{X_K/K} = (u_{X/k})_K.\]
The following proposition on base change of the Albanese torsor is known for separable extensions \cite[Corollary A.5]{wittenbergAlbaneseTorsorsElementary2008}, but seems not to have been recorded for inseparable extensions.
\begin{prop}\label{prop:albanesebasechange}
	Let $X$ be a geometrically integral variety over a field $k$ and $K/k$ a field extension. The canonical morphism
	\[ \theta^n_{K/k} \colon \Alb^n_{X_{K}/K} \longrightarrow \Alb^n_{X/k}\times_k K \]
	is a universal homeomorphism for $n=0,1$. Moreover, if $K/k$ is a (not necessarily algebraic) separable extension, then $\theta^n_{K/k}$ is an isomorphism for $n=0,1$.
\end{prop}
\bp 
Assume that $K/k$ is separable. Then the assertion is exactly \cite[Corollary A.5]{wittenbergAlbaneseTorsorsElementary2008}. Therefore, we may assume that $k$ is of characteristic $p>0$ and $K/k$ is purely inseparable.

Step 1: Let $K/k$ be a finite purely inseparable extension, then $\theta^n_{K/k}$ is surjective.

Since $X/k$ is geometrically integral, the smooth locus is dense and we find a finite separable extension $k_0/k$ such that we have a point $x \in X(k_0)$. Let $K_0=K k_0$ be the compositum in some algebraic closure of $k$. By the universal property of the Albanese, we have
\[  (\theta^n_{k_0/k})_{K_0}\circ \theta^n_{K_0/k_0}=      (\theta^n_{K/k})_{K_0}     \circ   \theta^n_{K_0/K}.\]
Since $k_0/k$ is separable, $\theta^n_{k_0/k}$ is an isomorphism. Therefore the surjectivity of $ \theta^n_{K_0/k_0}$ implies the surjectivity of $ (\theta^n_{K/k})_{K_0} $, which in turn implies the surjectivity of $\theta^n_{K/k}$. Thus after replacing $K/k$ with $K_0/k_0$, we may assume that $x\in X(k)$. Then $\Alb^1_{X/k}$ is a trivial $\Alb^0_{X/k}$-torsor and we obtain a maps 
\[ a\colon X \to \Alb^0_{X/k} \quad \text{and} \quad b \colon X_K \to  \Alb^0_{X_K/K} ,\]
both mapping $x$ to $0$ and satisfying
\[\theta^0_{K/k} \circ b = a_K.\]
In particular, surjectivity of $\theta^1_{K/k}$ is equivalent to the surjectivity of $\theta^0_{K/k}$. Let $B\subset \Alb^n_{X/k}\times_k K$ be the scheme-theoretic image of $\theta^0_{K/k}$. Since $\theta^0_{K/k}$ is a group homomorphism, $B$ is a subgroup.  Let $K \subset \Omega$ an algebraic closure. By \cite[Appendix A, Lemma A.2]{wittenbergAlbaneseTorsorsElementary2008} the map 
\[a_\Omega \colon X_{\Omega} \to  \Alb^0_{X/k}\times_k \Omega\]
is generating, which means that $\Alb^0_{X/k}\times_k \Omega$ is the smallest subgroup scheme over which $a_\Omega$ factors. Hence, we have 
\[B_\Omega = \Alb^0_{X/k}\times_k \Omega,\]
which implies $B= \Alb^n_{X/k}\times_k K$. We conclude that $\theta^0_{K/k}$ is surjective.

Step 2: The claim holds for a finite purely inseparable extension $K/k$.

Choose $m\geq 0$ such that $K^{p^m}\subset k$ and let
\[\phi \colon \spec(k) \longrightarrow\spec(K)\]
be the morphism induced by the $m$-th power of Frobenius. Let $A^n:=\phi^*\Alb^n_{X_K/K}$ be the pullback along $\phi$. We obtain a morphism
\[ \Phi \colon X \xrightarrow{F^m_{X/k}}X^{(p^m)} = \phi^*X_K \xrightarrow{\phi^*u_{X_K/K}} \phi^*\Alb^1_{X_K/K} =A^1.       \] The universal property of $\Alb^1_{X/k}$, implies the existence of unique $k$-morphisms 
\[f^n \colon \Alb^n_{X/k} \longrightarrow A^n, \qquad n=0,1 \]
such that $f^0$ is a group homomorphism, $f^1$ is $f^0$ equivariant and
\[f^1 \circ u_{X/k}=\Phi.\]
Let $F_n$ be the relative Frobenius of $\Alb^n_{X_K/K}$ over $K$. We have
\[ F^m_1\circ u_{X_K/K}=(u_{X_K/K})^{(p^m)}\circ F^m_{X_K/K}= \Phi_K = (f^1)_K \circ ( u_{X/k})_K= (f^1)_K \circ \theta^1_{K/k} \circ u_{X_K/K}, \]
which implies
\[ F^m_n = (f^n)_K \circ  \theta^n_{K/k}.\]
Since $F^m_n$ is a universal homeomorphism, $\theta^n_{K/k}$ is universally injective. By Step 1, we conclude that $\theta^n_{K/k}$ is bijective. The inverse map of underlying topological spaces is given by
\[       (\theta^n_{K/k})^{-1} =(F^m_n )^{-1}  \circ \theta^n_{K/k},   \]
which is continuous. Therefore $ \theta^n_{K/k}$ is a homeomorphism. The same argument applies to any base change, which proves that $\theta^n_{K/k}$ is a universal homeomorphism.

Step 3: The claim holds for any purely inseparable extension $K/k$. 

Since $K/k$ is algebraic, $K$ is the union of all finite subextensions $k'/k$ of $K/k$. Let $k' \subset k''$ be a subextension of $K/k$, which are finite over $k$. Consider the composition 
\[  \Alb^n_{X_{k''}/k''}\longrightarrow \Alb^n_{X_{k'}/k'}\times_{k'} k'' \longrightarrow\Alb^n_{X_{k'}/k'}.   \]
The first map is a universal homeomorphism by Step 1, hence affine and the second map is clearly affine. Therefore, the cofiltered limit
\[A^n:=\lim_{k'}\Alb^n_{X_{k'}/k'}    \]
exists as a $K$-scheme. Moreover, $A^0/K$ is a semi-abelian variety over $K$ and $A^1$ is a torsor under $A^0/K$. Since everything in sight is finitely presented, the universal property of the Albanese implies that the canonical map
\[ \Alb^n_{X_{K}/K} \lisim A^n\]
is an isomorphism. Therefore 
\[ \theta^n_{K/k} = \lim_i\theta^n_{k_i/k}\]
is a cofiltered limit of universal homeomorphisms, hence a universal homeomorphism.
\ep 
\begin{prop}[Spie\ss--Szamuely,  {\cite[Proposition 4.4]{spiessAlbaneseMapSmooth2003}}]\label{prop:albanesefinitekernel}
	Let $k$ be a field of characterstic $p\geq 0$ with separable closure $\bar{k}$. Let $\overline{X}/k$ be a smooth, projective, geometrically integral variety and $X\subset \overline{X}$ a dense open with Albanese torsor $X\to \mathrm{Alb}^1_{X/k}$. The induced homomorphism
	\[             \pi_1(X_{\bar{k}})^{\ab,(p')}\twoheadrightarrow \pi_1(\mathrm{Alb}^1_{X/k}\times_k\bar{k})^{(p')}   \]
	is surjective with finite kernel. 
\end{prop}
\bp 
Since $\bar{k}/k$ is separable,  \cref{prop:albanesebasechange} allows us to reduce to the case where $\bar{k}=k$. Let $K/k$ be an algebraic closure, we have a commutative diagram 
\[\begin{tikzcd}
	{\pi_1(X_K)^{\ab,(p')}} && {\pi_1(\Alb^1_{X_K/K})^{\ab,(p')}} \\
	&& {\pi_1(\Alb^1_{X/k}\times_kK)^{\ab,(p')}} \\
	{\pi_1(X)^{\ab,(p')}} && {\pi_1(\Alb^1_{X/k})^{\ab,(p')}.}
	\arrow["{u_{X_K/K}}", from=1-1, to=1-3]
	\arrow[from=1-1, to=3-1]
	\arrow["{\theta_{K/k}}", from=1-3, to=2-3]
	\arrow[from=2-3, to=3-3]
	\arrow["{u_{X/k}}", from=3-1, to=3-3]
\end{tikzcd}\]
By \cref{prop:albanesebasechange}, the map $\theta_{K/k}$ induces an isomorphism on fundamental groups. Since $K/k$ is purely inseparable, the homomorphisms 
\[ \pi_1(X_K)\lisim  \pi_1(X) \quad \text{and} \quad \pi_1(\Alb^1_{X/k}\times_kK) \lisim \pi_1(\Alb^1_{X/k}) \]
are isomorphisms. Hence, all vertical maps are isomorphisms. Therefore, after replacing $X/k$ with $X_K/K$, we may assume that $k$ is algebraically closed. Then the assertion follows from  \cite[Proposition 4.4]{spiessAlbaneseMapSmooth2003}.
\ep 
\subsection{Geometrically unibranch varieties}
We first prove an auxiliary lemma.
\begin{lem}\label{lem:induction}
	Let $k$ be a finitely generated field of characterstic $p\geq 0$ and $K/k$  a finite extension with separable closure $\overline{K}$. Let $\overline{X}/K$ be a smooth, projective, geomerically irreducible variety and $X\subset \overline{X}$ a dense open. Let $L/k$ be a Galois extension such that $\Gal(L/k)^\ab$ has finite exponent. Then 
	\[\left[ \mathrm{Ind}_{\Gal_{K}}^{\Gal_k}\pi_1(X_{\overline{K}})^{\ab,(p')}\right]_{\Gal_L}\]
	is finite.
\end{lem}
\bp 
Let $X\to \mathrm{Alb}^1_{X/K}=:A^1$ be the Albanese torsor of $X/K$. By \cref{prop:albanesefinitekernel}, we have an exact sequence of $\Gal_{K}$-modules
\[0 \to C\to \pi_1(X_{\overline{K}})^{\ab,(p')}\twoheadrightarrow \pi_1(A^1_{\overline{K}})^{\ab,(p')}\to 0, \]
with $C$ finite. Since induction is exact, we obtain an exact sequence 
\[0 \to  \mathrm{Ind}_{\Gal_{K}}^{\Gal_k}C\to  \mathrm{Ind}_{\Gal_{K}}^{\Gal_k}\pi_1(X_{\overline{K}})^{\ab,(p')}\twoheadrightarrow  \mathrm{Ind}_{\Gal_{K}}^{\Gal_k}\pi_1(A^1_{\overline{K}})^{\ab,(p')}\to 0. \]
Because $C$ is finite and applying $\Gal_L$-coinvariants is right-exact, we may replace $X$ with $A^1$ and assume that $X$ is a torsor under a semi-abelian variety $B$. In that case we have an isomorphism of $\Gal_{K}$-modules 
\[ \pi_1(X_{\overline{K}})^{(p')}= \pi_1(B_{\overline{K}})^{(p')}=\prod_{\l \neq p} T_\l B.\]
Therefore, it suffices to treat the case where $X=B$ is a semi-abelian variety over $K$.
Let $K/k_s/k$ be the maximal subextension such that $k_s/k$ is separable. Then $K/k_s$ is purely inseparable. Hence, there is an $m\geq 0$ such that $K^{p^m}\subset k_s$. Let $B_0/k_s$ be the semi-abelian variety obtained by base-changing $B
/K$ along 
\[K \to k_s, \quad x \longmapsto x^{p^m}.\] 
The morphism $B_0 \to B$, induces an isomorphism of $\Gal_{k_s}=\Gal_{K}$-modules
\[ T_\l B_0 \lisim T_\l B\]
for all $\l\neq p$. Thus, we have reduced to the case where $K=k_s$, i.e.\ $K/k$ is separable. When $K/k$ is separable, we have an isomorphism of $\Gal_k$-modules
\[\mathrm{Ind}_{\Gal_{K}}^{\Gal_k}T_\l B=T_\l \mathrm{Res}_{K/k}B,\]
where $\mathrm{Res}_{K/k}B$ is the Weil restriction. After replacing $B$ with its Weil restriction, we may further assume that $k=K$. The Tate-module of the semi-abelian variety $B/k$ lies in an exact sequence 
\[0 \to \prod_{\l \neq p}T_\l\mathbb{T}\to \prod_{\l \neq p}T_\l B \to \prod_{\l \neq p}T_\l A \to 0,\]
where $A/k$ is an abelian variety and $\mathbb{T}/k$ is a torus. Applying coinvariants is right exact, hence it suffices to treat the case where $B$ is a torus or an abelian variety. In both cases \cref{prop:coinvariantsofpure} applies to $\Lambda_\l:=T_\l B$. \ep 
\begin{thm}\label{thm:fabnessoffundamentalgroups}
	Let $X$ be a geometrically connected and geometrically unibranch variety over a finitely generated field $k$ of characteristic $p\geq 0$. Let $L/k$ be a (possibly infinite) Galois extension such that $\Gal(L/k)^\ab$ has finite exponent.
	\be 
	\item The group of coinvariants 
	\[\left[\pi_1(X_{\bar{k}})^{\ab,(p')}\right]_{\Gal_L}\]
	is finite.
	\item Let $\bar{x} \in X_{\bar{k}}$ be a geometric point and $s \colon \Gal_k\to \pi_1(X,\bar{x})$ a section of the exact sequence 
	\[1 \to \pi_1(X_{\bar{k}},\bar{x})\to \pi_1(X,\bar{x})\to \Gal_k\to 1.\] Endow $\pi_1(X_{\bar{k}},\bar{x})$ with a $\Gal_k$-action via conjugation with $s$. For any open subgroup $H \subset \left[\pi_1(X_{\bar{k}},\bar{x})\right]_{\Gal_L}$ the quotient  
	\[H^{\ab,(p')}\]
	is finite. In particular, if $p=0$, the group $\ \left[\pi_1(X_{\bar{k}},\bar{x})\right]_{\Gal_L}$ is FAb.
	\ee 
\end{thm}
\bp 
(1) We may replace $X$ with its associated reduced scheme $X_{\mathrm{red}}$, as this does not change the \'etale fundamental group and hence assume that $X$ is reduced. Since $X$ is geometrically unibranch, any non-empty open $U\subset X$ induces a surjective map
\[\pi_1(U_{\bar{k}}) \twoheadrightarrow\pi_1(X_{\bar{k}}) \]
and hence a surjective map on coinvariants. Thus, after shrinking $X$, we may additionally assume that $X/k$ is separated and integral. By \cite[Theorem 4.1 and Remark 4.2]{dejongSmoothnessSemistabilityAlterations1996}, we find an alteration $Y\to X$ such that $Y$ is a dense open in a projective variety $\overline{Y}/k$  and a finite extension $K/k$ such that $\overline{Y}$ is geometrically irreducible and smooth over $K$. Let $\overline{K}$ be the separable closure of $K$. Since $Y\to X$ is an alteration and $X$ is geometrically unibranch, the map 
\[\pi_1(Y_{\overline{K}})^{\ab,(p')} \longrightarrow\pi_1(X_{\overline{K}})^{\ab,(p')}=\pi_1(X_{\bar{k}})^{\ab,(p')}\]
is $\Gal_{K}$-equivariant and open. Thus, we obtain a $\Gal_k$-equivariant, open map 
\[ M:=\mathrm{Ind}_{\Gal_{K}}^{\Gal_k}\pi_1(Y_{\overline{K}})^{\ab,(p')} \longrightarrow  \pi_1(X_{\bar{k}})^{\ab,(p')}.\]
In particular, the induced map 
\[M_{\Gal_L} \longrightarrow  \left[\pi_1(X_{\bar{k}})^{\ab,(p')}\right]_{\Gal_L}\]
is open. The group $M_{\Gal_L}$ is finite by \cref{lem:induction}, hence $\left[\pi_1(X_{\bar{k}})^{\ab,(p')}\right]_{\Gal_L}$ is finite. 

(2) Let $\Pi:= \pi_1(X_{\bar{k}},\bar{x})$ and write 
\[\phi\colon \Pi \twoheadrightarrow \Pi_{\Gal_L}\]
for the quotient map. Let $H\subset  \Pi_{\Gal_L}$  be an open subgroup. Then $\phi^{-1}(H)$ is open in $\Pi$. We find an open normal subgroup $N\subset \pi_1(X)$ such that $\Pi':=N \cap   \Pi$ is contained in $\phi^{-1}(H)$. In particular, $\Pi'$ is stable under conjugation by $s$. Thus 
\[ \Pi' \rtimes s(\Gal_k) \subset \pi_1(X)\]
is an open subgroup, corresponding to a finite \'etale cover $Y\to X$ with 
\[ s(\Gal_k)\subset \pi_1(Y) \quad \text{and} \quad \pi_1(Y_{\bar{k}})= \Pi'.\] 
We let $\Gal_k$ act on $\Pi'=\pi_1(Y_{\bar{k}})$ via conjugation with $s$. The map
\[\pi_1(Y_{\bar{k}}) = \Pi' \subset \phi^{-1}(H) \twoheadrightarrow H^{\ab,(p')} \]
has  open image and factors through the quotient
\[\pi_1(Y_{\bar{k}}) \twoheadrightarrow \left[\pi_1(Y_{\bar{k}})^{\ab,(p')}\right]_{\Gal_L}, \]
which is finite by (1). Therefore $H^{\ab,(p')}$ is finite.
\ep 

\begin{rem}
 The assumption that $X/k$ is geometrically unibranch in \cref{thm:fabnessoffundamentalgroups} cannot be removed. Indeed if $X/k$ is the split nodal cubic, then $\Gal_k$ acts trivially on $\pi_1(X_{\bar{k}})\iso \widehat{\Z}$.
\end{rem}

\subsection{Arithmetic schemes}
Let $F$ be a number field and let $F^{\nr}/F$ be the maximal extension unramified at all finite places, so that $\pi_1(\o_F)=\Gal(F^{\nr}/F)$. By class field theory, $\pi_1(\o_F)^\ab=\Gal(F^{\nr}/F)^\ab$ is the narrow class group, which is finite, in particular has finite exponent. 

Let $\X$ be a normal scheme, flat, surjective
and of finite type over $\o_F$, with
geometrically connected generic fiber $\X_F$. We have a commutative diagram of fundamental groups 
\[\begin{tikzcd}
	1 & {\pi_1(\X_{\overline{F}})} & {\pi_1(\X_F)} & {\Gal_F} & 1 \\
	& {\pi_1(\X_{\overline{F}})} & {\pi_1(\X)} & {\pi_1(\o_F)} & 1.
	\arrow[from=1-1, to=1-2]
	\arrow[from=1-2, to=1-3]
	\arrow[no head, from=1-2, to=2-2]
	\arrow[shift left, no head, from=1-2, to=2-2]
	\arrow[from=1-3, to=1-4]
	\arrow[two heads, from=1-3, to=2-3]
	\arrow[from=1-4, to=1-5]
	\arrow[two heads, from=1-4, to=2-4]
	\arrow[from=2-2, to=2-3]
	\arrow[from=2-3, to=2-4]
	\arrow[from=2-4, to=2-5]
\end{tikzcd}\label{diagram}\]
The top row is exact and the vertical arrows are surjective, since $\X$ is normal. The complex of groups 
\[ \pi_1( \X_{\overline{F}})\to\pi_1(\X) \to \pi_1(\o_F)\to 1\]
is not exact in general, hence we consider the following two subgroups 
\[ \pi_1^{\geo}(\X/\o_F):= \im(\pi_1(\X_{\overline{F}})\to \pi_1(\X) )\subset \ker(\pi_1(\X)\to \pi_1(\o_F))=:K(\X/\o_F).\]
The subgroup $\pi_1^{\geo}(\X/\o_F)$ is normal in $ \pi_1(\X) $ and hence in $K(\X/\o_F)$.
\begin{lem}\label{lem:sectionimpliesexactness}
	Suppose $\X/\o_F$ admits a section. Then the sequence 
	\[ \pi_1( \X_{\overline{F}})\to\pi_1(\X) \to \pi_1(\o_F)\to 1\]
	is exact. In other words, we have $ \pi_1^{\geo}(\X/\o_F)=K(\X/\o_F)$.
\end{lem}
\bp 
A section $x \in \X(\o_F)$ yields compatible sections of $\pi_1(\X_F)\twoheadrightarrow \Gal_F$ and $\pi_1(\X)\twoheadrightarrow \pi_1(\o_F)$. The assertion then follows by a diagram chase as in \cite[Lemma 2.2]{luFinitenessEtaleFundamental2026}, whose proof does not require that $\X/\o_F$ is separated.
\ep 
We recall the following theorem of Rumely and Moret-Bailly \cite[Th\'eor\`eme 0.7.]{moret-baillyProblemesSkolemChamps2001}.
\begin{thm}[Rumely, Moret-Bailly]\label{thm:rumely}
	There is a finite extension $E/F$ such that $\X(\o_{E}) \neq \emptyset$.
\end{thm}
We also recall some well-known properties about FAb profinite groups.
\begin{lem}\label{lem:fabproperties}
	Let $G,\;K$ and $Q$ be profinite groups.
	\be 
	\item Let $K\to G$ be a homomorphism with open image. If $K$ is FAb, then $G$ is FAb.
	\item Let $K \subset G$ be an open subgroup. Then $K$ is FAb if and only if $G$ is FAb.
	\item Let $1\to K\to G \to Q\to 1$ be an extension. If $K$ and $Q$ are FAb, then so is $G$.
	\ee 
\end{lem}
\bp 
(1) A quotient of a FAb group is again FAb, hence (1) follows from (2).

(2) Let $K\subset G$ be an open subgroup. Assume that $K$ is FAb. Let $H\subset G$ be an open subgroup. Then $H':=H\cap K$ is open in $K$ and $H$, hence $(H')^\ab$ is finite and 
\[ (H')^\ab \to H^\ab\]
has open image. Therefore the image of the finite group $(H')^\ab$ has finite index in $H^\ab$. We conclude that $H^\ab$ is finite, hence $G$ is FAb. The converse is clear. 

(3) Let $1 \to K \xrightarrow{\alpha} G  \xrightarrow{\beta} Q \to 1$ be an exact sequence with $K$ and $Q$ FAb. Let $H\subset G$ be an open subgroup. Then we have an exact sequence 
\[1 \to H \cap \alpha(K) \to H \to \beta(H) \to 1. \]
The subgroups $H\cap \alpha(K) \subset \alpha(K) \iso K$ and $\beta(H)\subset Q$ are open, hence have finite abelianization. The sequence 
\[\left(H\cap \alpha(K)\right)^\ab \to H^\ab \to \beta(H)^\ab \to 1 \]
is exact, therefore $H^\ab$ is finite. Hence $G$ is FAb.
\ep 
\begin{thm}\label{thm:arithmeticfab}
	The group  $\pi_1^{\geo}(\X/\o_F)$ is FAb.
\end{thm}
\bp 
By \cref{thm:rumely}, there is a finite extension $E/F$ such that there is a section $x \in \X(\o_E)$. The induced map 
\[\pi_1^{\geo}(\X_{\o_E}/\o_E)\twoheadrightarrow\pi_1^{\geo}(\X/\o_F) \]
is surjective. \cref{lem:fabproperties} (1) implies that, after replacing $\X/\o_F$ with $\X_{\o_E}/\o_E$, we may assume $x \in \X(\o_F)$. The section $x$ induces compatible section of $\pi_1(\X_F,x)\twoheadrightarrow \Gal_F$ and $\pi_1(\X,x)\twoheadrightarrow\pi_1(\o_F)$.  Therefore the quotient map
\[\pi_1( \X_{\overline{F}},x)=\pi_1( \X_{\overline{F}})\twoheadrightarrow \pi_1^{\geo}(\X/\o_F) \]
factors through the group of coinvariants $\left[\pi_1( \X_{\overline{F}},x)\right]_{\Gal_{F^{\nr}}}$, which is FAb by \cref{thm:fabnessoffundamentalgroups}. Therefore, $\pi_1^{\geo}(\X/\o_F) $ is a quotient of a FAb group and hence is itself FAb.
\ep 
We recover the following result of Schmidt \cite[Corollary 3.3]{schmidtTameCoveringsArithmetic2002}, originally proved using results of Katz--Lang \cite{katzFinitenessTheoremsGeometric1981} and ramification theory.
\begin{cor}[Schmidt]\label{cor:schmidt}
The fundamental group $\pi_1(\X)$ is FAb.
\end{cor}
\bp 
By \cref{thm:rumely}, there is a finite extension $E/F$ such that $\X(\o_E)\neq \emptyset$. The induced homomorphism 
\[ \pi_1(\X_{\o_{E}})\to \pi_1(\X) \]
is open.  \cref{lem:fabproperties} (1) implies that, after replacing $\X/\o_F$ with $\X_{\o_E}/\o_E$, we may assume that $\X(\o_F)\neq \emptyset$. Then \cref{lem:sectionimpliesexactness} implies that the sequence
\[1 \to \pi_1^{\geo}(\X/\o_F)\to \pi_1(\X) \to \pi_1(\o_F)\to 1\]
is exact.
The group $ \pi_1^{\geo}(\X/\o_F)$ is FAb by \cref{thm:arithmeticfab} and $\pi_1(\o_F)$ is FAb by class field theory. We conclude that $\pi_1(\X) $ is FAb using \cref{lem:fabproperties} (3).
\ep 
\subsection{Examples of regular proper arithmetic surfaces}
Let $F$ be a number field and $\X$ be a normal scheme, flat, surjective
and of finite type over $\o_F$, with
geometrically connected generic fiber $\X_F$ We call $\X/\o_F$ an \textbf{arithmetic surface} if $\X\to \spec(\o_F)$ is of relative dimension 1.  

Let $\X/\o_F$ be an arithmetic surface and $\Y\to \X$ a connected finite \'etale cover. Then there is a unique finite extension $E/F$ such that $\Y \to \spec(\o_F)$ factors as
\[\Y \to \spec(\o_E)\to \spec(\o_F)\]
and $\Y/\o_E$ is an arithmetic surface. In fact, $E$ is the constant field of $\Y_F$. 
\begin{lem}\label{lem:quadraticextensionsurfaces}
	Let $E/F$ be a quadratic extension such that $2$ is unramified in $E$. There are infinitely many regular proper arithmetic surfaces $\X/\o_F$ admitting a connected degree 2 finite \'etale cover $\mathcal{Y}\to \X$ such that
	\begin{itemize}
		\item the morphism $\Y\to \spec(\o_F)$ factors over $\spec(\o_E)$ and
		\item $\Y/\o_E$ is an arithmetic surface with $\mathcal{Y}(\o_E)\neq \emptyset$.
	\end{itemize} 
\end{lem}
\bp 
For $n\geq 2$ consider the smooth affine curve $Z/\Q$ given by the equation
\[y^2=x^{2^n}+1.\]
Clearly $(0,1) \in Z(\Q)$. Let $Z \subset \overline{Z}$ be the smooth completion. Another affine chart of $\overline{Z}$ is given by 
\[v^2= u^{2^n}+1 \quad \text{with $x=1/u$ and $y=v/u^{2^{n-1}}$.} \] 
The minimal regular proper model $\mathcal{Z}$ of $\overline{Z}$ is an arithmetic surface over $\Z$ admitting a section. Define 
\[\mathcal{Y}:= \mathcal{Z} \times_{\Z} \spec(\o_E).\]
Then $\mathcal{Y}/\o_E$ is arithmetic surface with  $\Y(\o_E)\neq \emptyset$. The involution 
\[\tau\colon Z^\circ \longrightarrow Z^\circ,\quad  (x,y)\longmapsto (-x,-y)\] 
extends to an involution $\tau \in \Aut(\mathcal{Y}/\o_E)$, which has no fixed points on $\mathcal{Y}[1/2]$. Let $\sigma \in \Gal(E/F)$ be the unique non-trivial element and consider the automorphism 
\[\phi:=(\tau,\sigma) \in \Aut(\mathcal{Y}/\o_F).\]  Let $\p \in \spec(\o_F)$ be a prime above $2$. The automorphism $\sigma $ acts freely on $\spec(\o_E/\p)$, since $2$ and therefore $\p$ is unramified in $E$. Hence $\phi$ acts freely on $\mathcal{Y}$, and therefore the quotient map
\[\mathcal{Y}\to  \mathcal{Y}/\phi=:\X\]
is finite \'etale of degree 2. The geometric generic fiber 
 \[\X_{\overline{F}}\iso \overline{Z}_{\overline{F}}\]
 is connected, hence $\X/\o_F$ has the desired properties. For varying $n$, this construction yields distinct arithmetic surfaces, as the genus of the geometric generic fiber is $2^{n-1}-1$.
\ep 
For a profinite group $G$, we set 
\[d_2(G):= \dim_{\F_2} \H^1(G,\F_2).\] For a set of primes $S$ of a number field $F$, we write $\Gal_{F,S}$ for the Galois group of the maximal extension $F_S/F$ unramified outside $S$. Moreover $S_2(F)$ denotes the set of primes of $F$ lying above $2$.
\begin{lem}\label{lem:relativeinertianotfab}
	Let $F/\Q$ be a totally imaginary number field unramified at $2$ and $\widetilde{F}/F$ an infinite everywhere unramified pro-$2$ extension. Assume there exists a finite set of finite primes $T$ of $F$ splitting completely in $\widetilde{F}$ such that
	\[T \cap S_2(F)= \emptyset \quad \text{and} \quad \# T >d_2(\Cl(F)) +2[F\colon \Q] -1.\] 
	Let $E/F$ be a quadratic extension unramified at $2$ and ramified at every prime in $T$. Set
	\[I_{E/F}:= \ker(\pi_1(\o_E) \to \pi_1(\o_F))=\Gal(E^{\nr}/EF^{\nr}).\]
	Then the pro-$2$ completion $I_{E/F}^{(2)}$ is not finitely generated.
\end{lem}
\bp 
It suffices to prove that $d_2(I_{E/F})$ is infinite.

Step 1: Choose subextensions $\widetilde{F}/F_i/F$ such that $F_i/F$ is finite and Galois and the degree $n_i:=[F_i:F]\to \infty$ goes to infinity. Let $T_i$  denote the primes above $T$ in $F_i$. Since $T$ is completely split, we have  $\# T_i=n_i \# T $.  Define 
\[ \H_{i}:=\H^1(\Gal_{F_i,T_i},\F_2) \subset  \H^1(\Gal_{F_i,T_i \cup S_2(F) },\F_2)=:\H_{i,2} \subset F_i^\times/F_i^{\times,2}\]
and consider the map
\[ v_i \colon \H_{i,2} \longrightarrow \F_2^{T_i}, \quad b \longmapsto (v_\p(b) \bmod 2)_{\p \in T_i}. \]
For $b \in \H_{i}$ we have 
\[v_i(b) = 0 \iff  F_i(\sqrt{b})  \subset F^{\nr}.\]
Since $E/F$ is ramified at every prime in $T$, Abhyankhar's Lemma implies that for any $b \in \H_i$, we have $EF_i(\sqrt{b}) \subset E^{\nr}$. 
 Therefore we obtain a map
\[\phi_i \colon \H_i \longrightarrow \H^1(I_{E/F} ,\F_2), \quad  b \longmapsto (\sigma \mapsto \restr{\sigma}{EF_i(\sqrt{b})F^{\nr}}), \]
which factors as 
\[\tilde{\phi}_i\colon v_i(\H_i) \longrightarrow\H^1(I_{E/F} ,\F_2).\]
Suppose $0 \neq v_i(b) \in \ker(\tilde{\phi}_i)$. Then 
\[F_i(\sqrt{b}) \subset EF^{\nr}.\] Since $F_i(\sqrt{b})/F$ is ramified, we have
\[ F_i(\sqrt{b})F^{\nr} =EF^{\nr}.\]
For any $\p \in T_i$ the extension $EF^{\nr}/F^{\nr}$ is ramified at any prime over $\p$. Hence $F_i(\sqrt{b})$ must be ramified at $\p$. Therefore 
\[ v_i(b)= (1,\dots,1).\]
We conclude that $\dim_{\F_2} \ker(\tilde{\phi}_i)\leq 1$ and hence
\[ d_2(I_{E/F})=\dim_{\F_2} \H^1(I_{E/F},\F_2)\geq \dim_{\F_2}\im(\tilde{\phi}_i)\geq  \dim_{\F_2} v_i(\H_i)-1.\]
Therefore it suffices to prove that $\dim_{\F_2} v_i(\H_i)$ tends to infinity.

Step 2: We have
\[ \dim_{\F_2} v_i(\H_{i})  \geq \dim_{\F_2}  \im(v_i) - \dim_{\F_2} \H_{i,2}/ \H_{i}.  \]
Consider the map 
\[\mathrm{cl} \colon \F_2^{T_i}\longrightarrow \left( \Cl(F_i)/\langle [\mathfrak{q}], \mathfrak{q} \in S_2(F)\rangle \right)\otimes \F_2, \]
which for $\p \in T_i$ sends the corresponding basis vector to its ideal class $[\p ]$. The sequence 
\[ \H_{i,2}\xrightarrow{v_i} \F_2^{T_i} \xrightarrow{\mathrm{cl}}\left( \Cl(F_i)/\langle [\mathfrak{q}], \mathfrak{q} \in S_2(F)\rangle \right)\otimes \F_2,\]
is exact in the middle, therefore we have
\[ \dim_{\F_2}\im(v_i)\geq \# T_i - d_2(\Cl(F_i))= n_i \# T - d_2(\Cl(F_i)).\]
As $F_i/F$ is an everywhere unramified 2-extension, Schreier's inequality applies and we obtain
\[d_2( \Cl(F_i))-1 \leq n_i(d_2( \Cl(F))-1). \]
This provides the lower bound
\[ \dim \im(v_i) \geq |T_i| -  d_2(\Cl(F_i))\geq n_i( \# T-  d_2(\Cl(F)) +1)-1 . \]
Moreover, we have an injection
\[   \H_{i,2}/\H_i \hookrightarrow \prod_{\p  \in S_2(F_i)}\H^1(F_{i,\p},\F_2)/\H^1_{\nr}(F_{i,\p},\F_2),   \]
which implies  
\[ \dim_{\F_2}  \H_{i,2}/\H_i  \leq \sum_{\p  \in S_2(F_i)}^{}\left([F_{i,\p}:\Q_2]+1\right) \leq [F_i\colon \Q] +\# S_2(F_i) \leq 2[F_i:\Q] =2[F\colon \Q]n_i.\]
Altogether we obtain
\[ \dim v_i(\H_i) \geq n_i( \# T-d_2( \Cl(F)) -2[F:\Q]+1) -1 \to \infty. \qedhere\]
\ep 
\begin{thm}\label{thm:infinitekernel}
	Let $F$ be an imaginary quadratic field with discriminant $\Delta_F=- \l_1\cdot \dots \cdot \l_r$, where $r\geq 6$ and $\l_1,\dots,\l_r$ are distinct odd rational primes. There exist infinitely many regular proper arithmetic surfaces $\X/\o_F$ for which the pro-$2$ completion $K(\X/\o_F)^{(2)}$ is not finitely generated.
\end{thm}
\bp 
 Since $2 \nmid \Delta_F$, we conclude that $2$ in unramified in $F$. Because $F$ is totally imaginary, $\pi_1(\o_F)$ coincides with the Galois group of the maximal extension $F_{\emptyset}/F$ which is unramified everywhere. In particular, 
\[\pi_1(\o_F)^\ab=\Cl(F)\]
is the full class group. By Gauss's genus theory, we have
\[d_2(\pi_1(\o_F))= d_2(\Cl(F))= r-1\geq 5 > 2+2 \sqrt{2}=:\alpha_{F,\emptyset}.  \]
By \cite[Theorem 4]{hajirCuttingTowersNumber2021} there is an infinite everywhere unramified pro-2 extension $\widetilde{F}/F$ and an infinite set of primes $S$ of $F$ which completely split in $\widetilde{F}$. After shrinking $S$, we may further assume that $S$ does not contain any prime above $2$ and $\infty$. We now choose a large enough finite subset $T\subset S$ such that
	\[\# T >d_2(\Cl(F)) +2[F\colon \Q] -1.\] 
Using \cite[(9.2.8) Theorem]{neukirchCohomologyNumberFields2008}, we find a quadratic finite extension $E/F$ ramified in $T$ and unramified in $S_2(F)$. \cref*{lem:relativeinertianotfab} implies that the pro-$2$ completion of the group  
\[I_{E/F}:=\ker(\pi_1(\o_F)\to \pi_1(\o_E))\]
 is not finitely generated. By \cref{lem:quadraticextensionsurfaces}, we find infinitely many regular proper arithmetic surfaces $\X/\o_F$ admitting a degree $2$ finite \'etale cover $\Y\to \X$, such that the structure morphism 
\[\Y\to \spec(\o_F)\]
factors over $\spec(\o_E)$ and $\Y/\o_E$ is an arithmetic surface. We obtain a commutative diagram 
\[\begin{tikzcd}
	1 & {K(\Y/\o_E)} & {\pi_1(\Y)} & {\pi_1(\o_E)} & 1 \\
	1 & {K(\X/\o_F)} & {\pi_1(\X)} & {\pi_1(\o_F)} & 1
	\arrow[from=1-1, to=1-2]
	\arrow[from=1-2, to=1-3]
	\arrow[hook, from=1-2, to=2-2]
	\arrow[from=1-3, to=1-4]
	\arrow[hook, from=1-3, to=2-3]
	\arrow[from=1-4, to=1-5]
	\arrow[from=1-4, to=2-4]
	\arrow[from=2-1, to=2-2]
	\arrow[from=2-2, to=2-3]
	\arrow[from=2-3, to=2-4]
	\arrow[from=2-4, to=2-5]
\end{tikzcd}\]
with exact rows. Since $\Y \to \X$ is finite \'etale of degree $2$, the subgroup $\pi_1(\Y)\subset \pi_1(\X)$ is open and has index $2$. Therefore the subgroup
\[ G:=K(\X/\o_F)\cap \pi_1(\mathcal{Y}) \] is open in $K(\X/\o_F)$ and of index dividing $2$. In particular, $G$ is normal in $K(\X/\o_F)$ and the induced sequence 
\[1\to G^{(2)}\to K(\X/\o_F)^{(2)}\to K(\X/\o_F)/G \to 1\]
is exact by \cite[Lemma 2.7]{lieblichGeneratorsRelationsEtale2010}. Thus, it suffices to prove that $G$ has a non-finitely generated pro-$2$ quotient. From the diagram we obtain the exact sequence
\[1 \to K({\mathcal{Y}}/\o_E) \to  K(\X/\o_F) \cap \pi_1(\Y)  \to I_{E/F}\to 1, \]
which implies that the non-finitely generated group $I_{E/F}^{(2)}$ is a quotient of $ G$.
\ep

\bibliographystyle{amsalpha}
\bibliography{coinvariantsbib}
\end{document}